\documentclass{amsart} 
\usepackage{latexsym,amssymb,graphicx}

\newtheorem{thm}{Theorem}
\newtheorem{lem}[thm]{Lemma}
\newtheorem{prop}[thm]{Proposition}
 
\newtheorem{de}[thm]{Definition} 
\newtheorem{rem}[thm]{Remark}
\newtheorem{ex}[thm]{Example}

\newcommand{\BZ}{{\mathbb{Z}}}

\newcommand{\BO}{{\mathcal{O}}}
\newcommand{\BD}{{\mathcal{D}}}

\newcommand{\BJ}{{\mathcal{J}}}
\newcommand{\BI}{{\mathcal{I}}}

\newcommand{\BS}{{\mathcal{S}_p}}

\newcommand{\Vp}{{V_p}}

\newcommand{\Si}{{\Sigma}}

\DeclareMathOperator{\Int}{Int}

\DeclareMathOperator{\genus}{genus}

\DeclareMathOperator{\Lk}{Lk}
\begin{document}

\title{On the Frohman Kania-Bartoszynska ideal}

\author{ Patrick M. Gilmer}
\address{Department of Mathematics\\
Louisiana State University\\
Baton Rouge, LA 70803\\
USA}
\email{gilmer@math.lsu.edu}
\thanks{This research was partially supported by NSF-DMS-0203486}
\urladdr{www.math.lsu.edu/\textasciitilde gilmer/}

\date{August 18, 2005}

\begin{abstract} The Frohman Kania-Bartoszynska ideal is an invariant associated to a $3$-manifold with boundary and a prime $p \ge 5$. We give some estimates of this ideal. We also calculate this invariant  for some $3$-manifolds constructed by  doing surgery on a knot in the complement of another knot.
\end{abstract}

\maketitle

Let $p= 2d+1 \ge 5$ be a prime and let  
\begin{equation*}\label{cyclo}
\BO=\left\{
\begin{array}{cl}
\BZ[\zeta_p] &\text{if $p \equiv -1 \pmod{4}~,$} \\
\BZ[\zeta_{p},i]=\BZ[\zeta_{4p}] &\text{if $p \equiv 1 \pmod{4}~.$}\\
\end{array}\right.
\end{equation*}
Here and elsewhere, $\zeta_n$ denotes $e^{2 pi i/n}.$
If $M$ is an oriented connected closed 
$3$-manifold, define 
$I_p(M) = \BD \, \langle (M)\rangle_p$, where $\langle \ \rangle_p$ is the invariant defined in \cite{BHMV2},  $\BD=\langle S^3\rangle_p^{-1}$  and 
$M$ is given a $p_1$-structure with $\sigma$-invariant zero. Alternatively $M$ is given the weight zero \cite{G}. This 
is the same normalization as in \cite{MR}.
For example, $I_p(S^3)=1$. The invariant $I_p$ takes values in $\BO$ by \cite{Mu2,MR}.
Here is the definition of the Frohman Kania-Bartoszynska ideal invariant, as given in \cite{GM}. 

\begin{de}[\cite{FK}]  \label{def1} Given a connected $3$-manifold $N$ with boundary,
 let $\BJ_p(N)$  be the  ideal in $\BO$ generated by 
\[ 
\{ I_p(M)| \text{M is a closed connected oriented $3$-manifold containing $N$} \}. \]
In the case $p\equiv 1 \pmod{4},$ we  define
$\BJ_p^+(N)= \BJ_p(N)\cap \BZ[\zeta_p]$.
\end{de} 

 Frohman and Kania-Bartoszynska actually made this definition using the $SU(2)$ theory in place of the $SO(3)$ theory used here.  If we  replace $I_p$ in the above definition with $I_{2p}$ (in the sense of \cite{MR}) and replace $\BZ [\zeta_p]$ with $\BZ [\zeta_{8p}]$, then we obtain the ideal invariant of \cite{FK} exactly.  They gave some interesting examples of punctured $3$-manifolds that could not embed in other $3$-manifolds, even though this is not ruled out by homology or homotopy type.   We remark that results for punctured 
 $3$-manifolds follow immediately from the fact that quantum invariants of closed $3$-manifolds are multiplicative under connect sum and do not require the structure of a TQFT. They were able to show that the ideal invariant, associated to a certain Turaev-Viro  invariant,  was non-trivial for the union of two solid tori glued together by identifying neighborhoods of $(2,1)$ curves on their boundary.
This is a $3$-manifold with boundary a torus and first homology $\BZ \oplus \BZ_2.$   We also recover this estimate in Remark \ref{est}.

We note that if $p\equiv 1 \pmod{4},$ the ideal $\BJ_p(N)$ is generated by scalars, according to Definition \ref{def1},  which are either in $\BZ[\zeta_p]$ or in $i \BZ[\zeta_p]$.  Thus $\BJ_p(N)$ is generated over $\BO$ by $\BJ_p^+(N)$. 
Thus $\BJ_p^+(N)$ contains the same information as
$\BJ_p(N).$ We make the following definition so that 
the ideal we study lies in $\BZ[\zeta_p]$ in every case.
This is convenient as $\BZ[\zeta_p]$ has a simpler ideal structure  than $\BZ[i,\zeta_p]$. 
  
\begin{de} Let  $\BI_p(N)$ denote $\BJ_p(N)$ if $p\equiv -1 \pmod{4},$
and denote $\BJ_p^+(N)$ if $p\equiv 1 \pmod{4}.$
\end{de} 

The ideal $\BI_p(N)$ is hard to compute, in general,  because the definition involves the quantum invariants of infinitely many 
$3$-manifolds. 
One has the immediate:

\begin{prop} [\cite{FK}] \label{FK}  If $N_1$ embeds in $N_2$,  then $\BI_p(N_2) \subset \BI_p(N_1).$
\end{prop}

Let  $M^\circ$  denote $M$ punctured, i.e.  $M$ with an open 3-ball deleted.  If $M $ is closed,  $\BI_p( M^\circ)$ is the principal ideal generated by $ I_p(M).$   This follows from the fact the $I_p(M)$ is multiplicative under connected sum.

If $N$ has connected boundary,
 Masbaum and the author showed \cite{GM}, $\BJ_p(N)$ can be computed from the evaluation of a finite number of explicitly given skein elements  in $S^3$. Such evaluations can be algorithmically computed using \cite{MV}. Similar results could be obtained for  $N$ with disconnected boundary.  We do not explore this here.
 
 From now on,  we let $N$  denote a connected $3$-manifold, with a connected boundary denoted $\Si.$ We also let $h= 1-\zeta_p.$ We have that $h^{p-1}$ is a unit of $\BZ[\zeta_p]$ times $p.$  Let $\omega'$ denote the skein element 
 $h^{-(d-1)} \sum_{i=0}^{d-1} (-1)^i [i+1] e_i$ in the Kauffman Bracket skein module of $S^1 \times D^2.$ This is, up to multiplication by a unit from $\BO$,  the same as the surgery skein element $\omega$ of \cite{BHMV2}. Below we have modified \cite[16.6]{GM} slightly, replacing $\omega$ with $\omega'$, and $\BJ_p(N)$ with
$\BI_p(N)$. We also do not require that the complement of $H$ be a handlebody.
 
\begin{thm} [\cite{GM}] \label{algo} Suppose $N$ is presented as  surgery on a link $L$ in the complement of a handlebody $H$  in $S^3$. 
Pick a lollipop tree $G$ for $H$. $\BI_p(N)$
is the $\BZ[\zeta_p]$ ideal generated by  the evaluation of the
$\dim( \Vp(\Si))$
skein elements given by coloring $L$ by 
$\omega'$
and putting in a neighborhood of $G$ the basis elements described in \cite{GM} 
for $\BS(-\Si)$ indexed by small colorings of $G$.
\end{thm}

Since $\BZ[\zeta_p]$ is a Dedekind domain,  we have unique factorization of ideals. 

\begin{de} If $\mathfrak j$ is a nonzero ideal in $\BZ[\zeta_p]$, let 
$\nu_h(\mathfrak j)$ denote the exponent of $(h)$ appearing in the prime  factorization of $\mathfrak j$. So  
$\mathfrak j = (h)^{\nu_h(\mathfrak j)}\prod_{i=1}^{n} {\mathfrak p_i}^{e_i}$, where ${\mathfrak p_i}$ are prime ideals distinct from $(h)$.  In this case, let  $\breve {\mathfrak j}$ denote $\prod_{i=1}^{n} {\mathfrak p}^{e_i}$.  If  $\mathfrak j$ is zero, define $\nu_h(\mathfrak j)$ to be infinite, and $\breve{\mathfrak j}$ to be zero as well. 
\end{de}

Let $L_p= \BZ[\zeta_p,  \frac 1 p ]$. Using \cite[p. 180 4F]{Rib}, we have that $\breve {\mathfrak j}=  {\mathfrak j} L_p \cap  \BZ[\zeta_p].$ 
We can give a similar but somewhat easier formula for $\breve \BI_p(N)$, than given for $\BI_p(N)$ in Theorem \ref{algo}. This formula does not rely on  \cite{GM,GMW}. It only relies on the variant of the \cite{BHMV2}-theory over $K_p = \BO[ \frac 1 p ]$  discussed in \cite{G}.  Given  a not-necessarily-lollipop trivalent spine $G$ for a 
 handlebody $H$ with boundary $-\Si$, bases for the 
$K_p$-module $\Vp(-\Si) = \BS(-\Si) \otimes K_p$  are given in \cite{BHMV2} in terms of  $p$-admissible colorings of $G$.

\begin{thm}  \label{algo2} Suppose $N$ is presented as  surgery on a link $L$ in the complement of a handlebody $H$  in $S^3$.  Let $G$ be a trivalent spine of $H$.
Then $\breve \BI_p(N) = \breve {\mathfrak n}$ where $\mathfrak n$ is the $\BZ[\zeta_p]$ ideal  generated by 
the evaluations of the $\dim( \Vp(-\Si))$
skein elements given by coloring $L$ by 
$\sum_{i=0}^{d-1} (-1)^i [i+1] e_i$ 
and the $p$-admissible colorings of $G$ which comprise some basis  of  $\Vp(-\Si)$. \end{thm}

\begin{proof} $\BJ_p(N) K_p$ is  ${\mathfrak n}K_p$ by the same proof as \cite[16.6]{GM},  but using $\Vp(\Si)$ in place of $\BS(\Si),$ and
$\sum_{i=0}^{d-1} (-1)^i [i+1] e_i$ in place of $\omega$, as they are the same in $K_p$  up to multiplication by a unit of $K_p$. It then follows that $\BI_p(N) L_p$ is  ${\mathfrak n}L_p.$ Thus 
$\breve \BI_p(N) = \breve {\mathfrak n}$
\end{proof}

In \cite{GM}, $\BI_5$ of a family of homology circles (with torus boundaries) is computed using Theorem \ref{algo}. For an infinite subfamily of these, it follows that they cannot embed in $S^3$. For this non-embedding result, $\breve \BI_5$ , calculated using Theorem \ref{algo2}, suffices.
 We  consider below further examples calculated via Theorem \ref{algo} below. But first we discuss some more general results.

\begin{thm}\label{bettip} Let $\Si$ be the boundary of $N$. We have that:  \(  \dim( H_1(N, \BZ_p) ) \ge \genus(\Si).\)
If   $ \dim( H_1(N, \BZ_p))= \genus(\Si)$, then  $\nu_h( \BI_p(N)=0$.
If  $ \dim( H_1(N, \BZ_p)) > \genus(\Si)$,
 $\nu_h( \BI_p(N) \ge d-1$.  Moreover
 \[ \nu_h( \BI_p(N) \ge  \frac {d-1} 3 ( \dim( H_1(N, \BZ_p) -\genus(\Si) ). \]     \end{thm}

\begin{proof} By duality,  the kernel of the map induced by inclusion: $j: H_1(\Si ,\BZ_p) \rightarrow H_1(N,\BZ_p)$ has dimension $\genus(\Si)$. See the argument in proof of \cite[Lemma 8.15]{L} for example. Thus the image of this map has image of dimension $\genus(\Si)$.

Now suppose that $\genus(\Si) = \dim( H_1(N, \BZ_p))$. Then we may pick a simple closed curve $\alpha$ which represents an element of the kernel of $j$ above. If we attach a handle to $N$ along  $\alpha$ to form a new $3$-manifold $N'$ with boundary $\Si'$, we find that 
$\genus(\Si') = \dim( H_1(N', \BZ_p)),$ but $\genus(\Si') =\genus(\Si)-1.$ Repeating this process, we reach a $\BZ_p$-homology ball  with boundary a 2-sphere. Capping this off with a 3-ball,  we obtain a  $\BZ_p$-homology sphere $M$ which contains $N$.
By \cite{Mu2},  $h$ does not divide $I_p(M)$, and so  $\nu_h( \BI_p(N))=0$.

Suppose $N$ embeds in $M$.  As the image of  the map induced by inclusion from $H_1(\Si, \BZ_p)$ to $H_1(N,Z_p)$ has dimension $\genus(\Si),$   using the Mayer-Vietoris sequence
\[ H_1(\Si, \BZ_p) \rightarrow H_1(N,Z_p) \oplus H_1(M \setminus \Int(N),Z_p) 
\rightarrow H_1(M,Z_p)  \rightarrow 0, \] we see  that 
\[ \dim (H_1(M,Z_p)) \ge  \dim( H_1(N, \BZ_p)) -\genus(\Si) .\]

If  $ \dim( H_1(N, \BZ_p)) > \genus(\Si)$, then any $3$-manifold $M$ containing $N$ has  \newline $H_1(M,Z_p) \ne 0$. By \cite[Thm 4.3]{CM}, $h^{d-1}$ divides $I_p(M)$.

The final result also follows in the same way from \cite[Thm 4.3]{CM}. \end{proof}

The cut number of $N$, denoted $c(N)$, is the maximal number of closed surfaces that one can remove from the interior of $N$ without disconnecting $N$.
Using  \cite[Theorem 15.1]{GM} we  also have,
\begin{thm}\label{cut}  \( \nu_h( \BI_p(N) \ge   (d-1) \  c(N). \)   \end{thm}

We wish to study $\BI_p$ of $3$-manifolds obtained by the following simple construction.  We start with  ordered link $L$ of two components $K$ and $J$ in  $S^3.$ We consider  $L_k$,  the result of surgery to the exterior of $J$ along $K$ with framing $k.$ 

\begin{lem}\label{lemma} We have that $H_1(L_k)  \approx \BZ \oplus \BZ_{\gcd(k,\Lk(J,K))}$, where  we interpret $\BZ_0$ as $\BZ$, and $\BZ_1$ as the trivial group. So $\dim (H_1(L_k,Z_p))>1$ if and only if 
\linebreak $\gcd(p,  k,\Lk(J,K)) \ne 1.$  Thus  $\nu_h( \BI_p(L_k) \ge d-1$ if and only if $\nu_h( \BI_p(L_k)>0$ if and only if   ${\gcd(k,\Lk(J,K),p) \ne 1}$. 
\end{lem}

\begin{proof} $H_1( S^3 \setminus J \setminus K)$ is free abelian on the meridians $m(J)$ and $m(K)$ of $J$ and $K$. The  filling along $K$ adds the relation $\Lk(J,K)\ m(J)+ k\ m(K)=0.$ Then we apply Theorem \ref{bettip}.
\end{proof}

In view of this lemma, we make the following: 

\begin{de} Let $L$ be a link with 2-components $J$ and $K$ as above.
If  $p$ does not divide $\gcd(k, \Lk(J,K)),$ and 
$\BI_p(L_k) =  (1)$, we  say $\BI_p(L_k)$ is large.
If  $p$  divides $\gcd(k, \Lk(J,K)),$ and 
$\BI_p(L_k) = (h^{d-1})$, we  also say $\BI_p(L_k)$ is large.
Otherwise we  say that $\BI_p(L_k)$ is small.
\end{de}

When $K$ is unknotted,  we do not need Theorem \ref{algo} to calculate $\BI_p(L_k)$ as the following theorem demonstrates.  
We remark that, if $J$ is unknotted, we may change $L$ by performing some full twists on the exterior of $J$. Let $L(s)$ denote the result of applying $s$ twists. Then $L_k$ is homeomorphic to $L(s)_{k +s(\Lk(J,K))}.$ 
Thus the following theorem also holds for links obtained from links where $J$ and $K$ are trivial by twisting around $J$ as described  above. 

\begin{thm} \label{large} If L is a link of two components with $K$ unknotted,  then $\BI_p(L_k)$ is large for all $p$.
\end{thm}

\begin{proof} Since $K$ is unknotted, $L_k$ embeds in the lens space $L(k,1)$ if $k \ne 0$ and in $S^1 \times S^2$ if $ k=0.$ 

To say  $p$ does not divide $\gcd(k, \Lk(J,K))$ is to say 
$p$ does not divide $k$ or to say $p$ does  divide $k$ but does not divide $\Lk(J,K).$

Suppose first that $p$ does not divide $k,$ then we have that
$I_p( L(k,1))$ is a unit of $\BZ[\zeta_p]$. This follows from \cite[Theoreme 1]{Ra} and \cite[Corollary 4.8]{MOO}. This implies that $\BI_p(L_k)= (1)$

Now suppose $p$ does divide $k$. As $S^1 \times S^2$ is obtained by 0-framed surgery on the unknot, and $L_p(L(k,1))$ is obtained by surgery on this same framed knot after inserting a multiple of $p$ twists, it follows that  $I_p(L(k,1))$ , $I_p(S^1 \times S^2)$ and $h^{d-1}$ all differ only by a factors which are units from $\BO.$ Thus $\BI_p(L_k)\supset  (h^{d-1})$. 

If we suppose that $p$ does not divide $\Lk(J,K),$ then  $\gcd(p,k,\Lk(J,K))=1,$ and so  by Lemma \ref{lemma}, $\nu_h(\BI_p(L_k))=0.$  This together with  $\BI_p(L_k)\supset  (h^{d-1})$ 
(which we know assuming that  $p$ does divide $k$) implies that $\BI_p(L_k)= (1).$

Now suppose that $p$ divides $\gcd(k, \Lk(J,K))$. As $p$ divides $k$,  $\BI_p(L_k)\supset  (h^{d-1})$.  As $\gcd(p, k, \Lk(J,K)) \ne 1$, by Lemma \ref{lemma}, $\BI_p(L_k)\subset  (h^{d-1}).$
  \end{proof}
  
  Let $K(k)$ denote $k$-framed surgery to $S^3$ along 
  $K$. By a similar proof, we have:
    
  \begin{thm} \label{large2} If $I_p(K(k))$ is a unit, then $p$ does not divide $k$ and  $\BI_p(L_k)$ is large. If $I_p(K(k))$ is up to units $h^{d-1}$ and $p$ divides $\Lk(J,K)$, then $p$ does divide $k$ and $\BI_p(L_k)$ is large.
  \end{thm}
  
  \begin{prop}\label{period} \(\BI_p(L_k) =\BI_p(L_{p+k}) \) \end{prop}

\begin{proof} Changing the framing of $K$ by $p$ leaves each term
in the expansion of a generator for the ideal the same since changing  the framing of a strand colored $i$ multiplies  an evaluation by a $p$th root of unity.
\end{proof}

By Theorem \ref{algo},
we have:
\begin{equation} \BI_p(L_k)=\left( \{ < K(t^k \omega') \cup J(v^s)> | 
0 \le s \le d-1\} \right), \label{videal} \end{equation}
where we use parenthesis to indicate the ideal generated by a set,  $K(t^k \omega)$
denotes $K$ cabled by the skein $t^k \omega$ and  $J(v^s)$  denotes $J$ cabled by $v^s.$ Here we use $t$ to denote the full twist map on a solid torus and the map this induces on the skein of the solid torus and $v$ denotes the skein element of the solid torus defined in \cite{GMW,GM}. This is the formula we  use below but we obtain another interesting formula if  we use instead the $\omega$ basis for $\BS$ of a torus \cite{BHMV2}: 
$ \{ t^i \omega | 0 \le s \le d-1\}$. 
\begin{equation} \BJ_p(L_k)=\left( \{ < K(t^k \omega) \cup J(t^s \omega)> | 
0 \le s \le d-1\} \right), \label{wideal} \end{equation}
Let $L_{k,s}$ denote the result of $k$-framed surgery along $K$  and $s$-framed surgery along $J$. Equation \ref{wideal}  says that $\BJ_p(L_k)$ is generated by 
\[ \{ I_p( L_{k,s}) |
0 \le s \le d-1\} \]
As the $L_{k,s}$ clearly include $L_k$, in this case,  $\BJ_p(L_k)$, the ideal generated by $I_p$ of all closed $3$-manifolds 
containing $L_k$,  is actually generated by $I_p$ of exactly $d$ of them.

We  used Bar-Natan's Knot Atlas \cite{B} to look for links $L$ in Thistlethwaite's list of prime links with low crossing number  with  $\BI_5(L_k)$ small for some $k$.  There is no  reason  to restrict oneself to considering prime links here, but it is convenient to use this list. We checked all links with two components and with up to  nine crossings.  In this range, all the links have at most one knotted component, and we took $K$ to be this component.  We also used Bar-Natan's mathematica package KnotTheory \cite{B} to help calculate $\BI_5(L_k).$ 

 \begin{figure}[ht]
\includegraphics[width=1in]{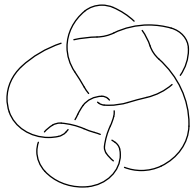}
\caption{L9a6 as drawn by Redelmeier's DrawPD command \cite{B} } \label{9a6}
\notag
\end{figure}

Many links 
were guaranteed to have large ideals by Theorem \ref{large} either because both components were unknotted or because one component was unknotted and by twisting about this unknotted component, one  could unknot the other component. There were many other other links that did not fit this criterion  but nevertheless for which we calculated that $\BI_5(L_k)$ were large using Theorem \ref{algo}, i.e. Equation \ref{videal}. The only  small 
$\BI_5(L_k)$ that we found  among the prime two component links with  nine or less crossings are listed in  following table.  The knots $K$ are denoted as in Rolfsen's table. The actual ambient orientation of $L$ is that encoded by the PD code in Bar-Natan's  package.  Linking number refers to the absolute value of the linking number of the two components of $L$.
Note that since $\BZ[\zeta_5]$ is a PID, $\BI_5(L_k)$ is necessarily principal. Small ideals appear to be rare among $3$-manifolds constructed in this way from prime links with low crossing number, but they do occur. The fourth line  describes the example given in \cite[16.7,16.8]{GM}.  Reading the first line for instance, we get  Example \ref{ex}. The other lines lead to similar examples. Similar data from a search through prime links with less than 12 crossings can be found at the author's homepage.

\begin{center}

\begin{tabular}{| l | c | c| c| c| c| c|}
\hline
L    & K       &  k & generator of   & norm of & linking & $L_k$ is a    \\ 
           &          &   &  $\BI_5(L_k)$             & generator & number&  homology $S^1$   \\ 
           \hline
L9a6   & $5_1$  & 5n  & $1- 2\ {\zeta_5}^2$ & 31 & 2 & if n odd   \\  \hline

L9a7   & $5_2$  & 5n+1  & $1 + 2\  {\zeta_5}^2$ & 11 & 2 &  if n even \\  \hline

L9a11   & $4_1$  & 5n+1  & $1 +2\  {\zeta_5}$ & 11 & 2 & if n even  \\  \hline

L9a12   & $5_1$  & 5n  & $1 + 2\  {\zeta_5}^2$ & 11 & 2 & if n odd  \\  \hline

L9a15   & $5_2$  & 5n+4  & $1 + 2\  {\zeta_5}^3$ & 11 & 0 & if n= -1  \\  \hline

L9a17   & $3_1$  & 5n+3  & $1 + 2\  {\zeta_5}^2$ & 11 & 0 & never  \\  \hline

L9a23   & $3_1$  & 5n+3  & $1 + 2\  {\zeta_5}^2$ & 11 & 3 &  if $n \ne 0 \mod{3}$  \\  \hline

L9a23   & $3_1$  & 5n+4  & $1 + 2\  {\zeta_5}$ & 11 & 3 &  if $n \ne 1 \mod{3}$  \\  \hline
\end{tabular}
\end{center}

\medskip

\begin{ex}\label{ex} Let $L$ be the alternating link pictured in  Figure \ref{9a6}. 
Let $K$ be the knotted component, and $J$ the unknotted component. $\BJ_5(L_k)$ is trivial unless $k \equiv 0\pmod{5}$. $\BJ_5(L_{5n})= (1 -2 {\zeta_5}^2).$ This is non-trivial as the norm of $1 -2 {\zeta_5}^2$ is 31. Thus $\BI_p(L_{5n})$ does not embed in any closed 
$3$-manifold $M$ (like the 3-sphere) for which $I_5(M)$ is not divisible by 
$1 -2 {\zeta_5}^2.$ If $n$ is odd,   $L_{5n}$ is a homology circle.
\end{ex}

We remark that finitely generated projective modules $S_{2p}(\Si)$ over $\BZ[\zeta_{8p}]$, analogous to $S_p(\Si)$, but associated to the SU(2) theory are defined in \cite[\S 5]{G}. If one could find an explicit finite set of generators for $S_{2p}(\Si)$, one could give an algorithm, as above,  to calculate the original  Frohman Kania-Bartoszynska ideal.  As \cite{BHMV2} describes a basis for  $V_{2p}(\Si)$ defined  over $\BZ[\zeta_{8p}, \frac 1{2p}]$, one can give 
 such a formula for the $\BZ[\zeta_{8p}, \frac 1 {2p}]$-ideal generated by $I_{2p}$ of all closed oriented $3$-manifolds containing $N$.  Subsequent to the first version of this paper, 
 Qazaqzeh \cite{Q} has  found a basis for $S_{2p}(\text{torus})$ and has used it to study  
the original  Frohman Kania-Bartoszynska ideal.

We can already calculate some of the ideals estimated/discussed in  \cite{FK}. 
Let $\tau_3(M)$ denote the quantum invariant considered by Kirby and Melvin
and denoted by them in this way.  It takes values in $\BZ[ \zeta_8]$ according to \cite[6.3]{KM}. Let $\BJ(N,\tau_3)$ denote the $\BZ[ \zeta_8]$-ideal generated by the set
$ 
\{ \tau_3 (M)|\text{M is a closed oriented }$ \linebreak $\text{$3$-manifold containing $N$} \} .$
Let $\BJ(N,TV_3)$ denote the  $\BZ[ \sqrt{2}]$-ideal generated by the set
$\{\tau_3 (M) \cdot \overline{\tau_3 (M)} |\text{M is a closed oriented  $3$-manifold containing $N$} \} .$
The ideals obtained from $\BJ(N,\tau_3)$ and $\BJ(N,TV_3)$ by extending rings to $\BZ[\zeta_{24}]$  are denoted by $I_3(N)$ and ${I_3}^{TV}(N)$ in \cite{FK}.

\begin{prop} \label{tau3} Let  L be a two component link. 
If  the linking number is odd, then
 $\BJ(L_{k},\tau_3)=  \BZ[ \zeta_8]$
and  $\BJ(L_{k},TV_3)= \BZ[ \sqrt {2}  ].$
These equalities also hold if  the linking number is even and $k$ is odd.
If the linking number is even and $k \equiv 0 \pmod{4}$,  then
 $\BJ(L_{k},\tau_3)= (\sqrt{2})  \subset \BZ[ \zeta_8]$
and 
$\BJ(L_{k},TV_3)= (2) \subset \BZ[ \sqrt {2}  ].$
If the linking number is even and $k \equiv 2 \pmod{4}$,  then
 $\BJ(L_{k},\tau_3)= (0)  \subset \BZ[ \zeta_8]$
and 
$\BJ(L_{k},TV_3)= (0) \subset \BZ[ \sqrt {2}  ].$  \end{prop}

 \begin{proof}   If the linking number is odd, consider $M$ obtained from $S^3$ by doing $0$-framed surgery along  $J$ and $k$-framed surgery along $K$.  By  \cite[6.1]{KM}  $\tau_3(M)=1$.   If the linking number is even, consider $M$ obtained by doing $1$-framed surgery on $J$ and the given surgery on $K$. If $k$ is odd, we again obtain $\tau_3(M)=1$. If  $k \equiv 0 \pmod{4}$, \cite[6.1]{KM}  we have that $\tau_3(M)=\sqrt {2}$, up to units.
Consult \cite[6.3]{KM} to see only possible values for any $M$ are multiples of $\sqrt {2}$.

For the last case,  any $3$-manifold containing $L_{k}$ can be obtained by $k$-framed surgery along $K$ and together with some integral surgery along a link contained in a tubular neighborhood of  $J$. These extra components must link $K$ an even number of times. As  $k \equiv 2 \pmod{4}$, the Brown invariant associated to the linking matrix for this surgery is undefined, and $\tau_3$ of the result is zero.
\end{proof}

\begin{rem}\label{est} {\em Let $T(4,2)$ denote the $(4,2)$ torus link. It consists of two unknots which link each other twice in a very simple  way.  We remark that $T(4,2)_0$  is homeomorphic to the union of two solid tori glued together by identifying neighborhoods of $(2,1)$ curves on their boundary.   Proposition  \ref{tau3} confirms the estimate in \cite{FK} that ${I_3}^{TV}(T(4,2)_0) \subset (2)$.}
\end{rem}

\end{document}